\documentclass[12pt,leqno,draft]{article}

\def\tr{\mathop{\rm tr}}

\newtheorem{theorem}{Theorem}
\newtheorem{lemma}[theorem]{Lemma}
\newtheorem{proposition}[theorem]{Proposition}
\newtheorem{definition}[theorem]{Definition}
\newtheorem{corollary}[theorem]{Corollary}

\newcommand{\begintheorem}{\addtocounter{equation}{1}\begin{theorem}}
\newcommand{\beginlemma}{\addtocounter{equation}{1}\begin{lemma}}
\newcommand{\beginproposition}{\addtocounter{equation}{1}\begin{proposition}}
\newcommand{\begindefinition}{\addtocounter{equation}{1}\begin{definition}}
\newcommand{\begincorollary}{\addtocounter{equation}{1}\begin{corollary}}

\begin{document}

\title{Some notes about matrices, 5}

\author{Stephen William Semmes 	\\
	Rice University		\\
	Houston, Texas}

\date{}

\maketitle

	As usual, let ${\bf R}$, ${\bf C}$ be the real and complex
numbers, and for each positive integer $n$ let ${\bf R}^n$, ${\bf
C}^n$ be the real and complex vector spaces of $n$-tuples of real and
complex numbers, respectively.  We write $GL({\bf R}^n)$, $GL({\bf
C}^n)$ denote the general linear groups of invertible real and
complex-linear transformations on ${\bf R}^n$, ${\bf C}^n$,
respectively, and $SL({\bf R}^n)$, $SL({\bf C}^n)$ for the subgroups
of $GL({\bf R}^n)$, $GL({\bf C}^n)$ consisting of linear
transformations with determinant equal to $1$.  Also, $O({\bf R}^n)$,
$U({\bf C}^n)$ are the subgroups of $GL({\bf R}^n)$, $GL({\bf C}^n)$
of orthogonal and unitary linear transformations on ${\bf R}^n$, ${\bf
C}^n$, respectively, which are the linear transformations which
preserve the standard Euclidean norm, or, equivalently, the standard
inner products on ${\bf R}^n$, ${\bf C}^n$, which is the same as
saying that they are invertible linear transformations whose inverses
are equal to their adjoints, and $SO({\bf R}^n)$, $SU({\bf R}^n)$ are
the special orthogonal and unitary groups, which are the subgroups of
$O({\bf R}^n)$, $U({\bf C}^n)$ of linear transformations which also
have determinant equal to $1$.

	Let us write $\mathcal{S}({\bf R}^n)$, $\mathcal{S}({\bf
C}^n)$ for the real vector spaces of self-adjoint linear
transformations on ${\bf R}^n$, ${\bf C}^n$.  We write
$\mathcal{S}_+({\bf R}^n)$, $\mathcal{S}_+({\bf C}^n)$ for the open
convex cones in $\mathcal{S}({\bf R}^n)$, $\mathcal{S}({\bf C}^n)$
consisting of linear transformations which are positive-definite.
Furthermore, we write $\mathcal{M}({\bf R}^n)$, $\mathcal{M}({\bf
C}^n)$ for the subsets of $\mathcal{S}_+({\bf R}^n)$,
$\mathcal{S}_+({\bf C}^n)$ of linear transformations with determinant
equal to $1$.

	Of course $GL({\bf R}^n)$, $GL({\bf C}^n)$ are open subsets of
the vector spaces of all linear transformations on ${\bf R}^n$, ${\bf
C}^n$, respectively, which can be defined by the condition that the
determinant is nonzero.  One can think of $SL({\bf R}^n)$, $SL({\bf
C}^n)$ as hypersurfaces in $GL({\bf R}^n)$, $GL({\bf C}^n)$, and
indeed they are regular or smooth hypersurfaces because the equation
that the determinant be equal to $1$ is nondegenerate on these
subsets, which is to say that the gradient of the determinant, as a
function on the vector space of linear transformations on ${\bf R}^n$
or ${\bf C}^n$, is nonzero at every invertible linear transformation.
Similarly, $\mathcal{M}({\bf R}^n)$, $\mathcal{M}({\bf C}^n)$ are
smooth hypersurfaces in the open sets $\mathcal{S}_+({\bf R}^n)$,
$\mathcal{S}_+({\bf C}^n)$.

	There are natural actions of $GL({\bf R}^n)$, $GL({\bf C}^n)$
on $\mathcal{S}({\bf R}^n)$, $\mathcal{S}({\bf C}^n)$, defined by
\begin{equation}
	A \mapsto T^* \, A \, T
\end{equation}
when $A$ is a self-adjoint linear transformation on ${\bf R}^n$ or
${\bf C}^n$ and $T$ is an invertible linear transformation on the same
space.  If $T$ is an invertible linear mapping on ${\bf R}^n$ or ${\bf
C}^n$, then this action takes $\mathcal{S}_+({\bf R}^n)$ onto itself
or $\mathcal{S}_+({\bf C}^n)$ onto itself, as appropriate.  Similarly,
if $T$ lies in $SL({\bf R}^n)$ or in $SL({\bf C}^n)$, then this action
takes $\mathcal{M}({\bf R}^n)$ onto itself or $\mathcal{M}({\bf C}^n)$
onto itself, as appropriate.

	The exponential of a linear transformation $A$ on ${\bf R}^n$
or on ${\bf C}^n$ is defined by
\begin{equation}
	\exp A = \sum_{j=0}^\infty \frac{A^n}{j!},
\end{equation}
and is always an invertible linear transformation on the same space.
When $A$ is self-adjoint, $\exp A$ is self-adjoint and
positive-definite, and indeed the exponential defines a one-to-one
mapping of $\mathcal{S}({\bf R}^n)$, $\mathcal{S}({\bf C}^n)$ onto
$\mathcal{S}_+({\bf R}^n)$, $\mathcal{S}_+({\bf C}^n)$, respectively.
If we write $\mathcal{S}_0({\bf R}^n)$, $\mathcal{S}_0({\bf C}^n)$ for
the linear subspaces of $\mathcal{S}({\bf R}^n)$, $\mathcal{S}({\bf
C}^n)$ consisting of self-adjoint linear transformations with trace
$0$, then the exponential is a one-to-one mapping of
$\mathcal{S}_0({\bf R}^n)$, $\mathcal{S}_0({\bf C}^n)$ onto
$\mathcal{M}({\bf R}^n)$, $\mathcal{M}({\bf C}^n)$, respectively.

	Next we look at some of these objects in terms of vector
calculus.

	For linear transformations on ${\bf R}^n$ or on ${\bf C}^n$,
the differential of the exponential mapping at the origin is equal to
the identity.  Explicitly, for each linear transformation $A$ on ${\bf
R}^n$ or on ${\bf C}^n$, we have that
\begin{equation}
	\biggl(\frac{d}{dt} \exp (t \, A) \biggr)_{t = 0} = A.
\end{equation}
Similarly, for the second differential, if $A$, $B$ are linear
transformations on ${\bf R}^n$ or on ${\bf C}^n$, then
\begin{equation}
	\biggl(\frac{\partial^2}{\partial s \partial t}
			\exp (s \, A + t \, B) \biggr)_{s, t = 0}
		= \frac{1}{2}(A \, B + B \, A).
\end{equation}

	Let us define a Riemannian metric on $\mathcal{S}_+({\bf
R}^n)$, $\mathcal{S}_+({\bf C}^n)$ as follows.  If $P$ is a
positive-definite self-adjoint linear transformation on ${\bf R}^n$ or
on ${\bf C}^n$, then tangent vectors to $\mathcal{S}_+({\bf R}^n)$,
$\mathcal{S}_+({\bf C}^n)$ at $P$ are given simply by self-adjoint
linear transformations $A$, $B$ on ${\bf R}^n$ or on ${\bf C}^n$, and
we define the inner product of these two tangent vectors at $P$ by
\begin{equation}
	\langle A, B \rangle_P = \tr (P^{-1} \, A \, P^{-1} \, B),
\end{equation}
where $\tr C$ denotes the trace of a linear transformation $C$.  Using
standard properties of the trace it is easy to see that this is a
positive-definite symmetric bilinear form in $A$, $B$, and it depends
smoothly on $P$, so that we get a Riemannian metric on
$\mathcal{S}_+({\bf R}^n)$ or on $\mathcal{S}_+({\bf C}^n)$, as
appropriate.

	Suppose that $T$ is an invertible linear mapping on ${\bf
R}^n$ or on ${\bf C}^n$, and consider the mapping $\tau = \tau^T$ on
$\mathcal{S}_+({\bf R}^n)$ or $\mathcal{S}_+({\bf C}^n)$, as
appropriate, defined by
\begin{equation}
	\tau(P) = T^* \, P \, T.
\end{equation}
If $A$, $B$ are self-adjoint linear transformations on ${\bf R}^n$ or
on ${\bf C}^n$, which we view as tangent vectors to
$\mathcal{S}_+({\bf R}^n)$ or $\mathcal{S}_+({\bf C}^n)$ at $P$, then
$d \tau_P(A)$, $d \tau_P(B)$, which are the images of $A$, $B$ as
tangent vectors at $P$ under the differential of the mapping $\tau$,
are self-adjoint linear transformations on ${\bf R}^n$ or ${\bf C}^n$
given by
\begin{equation}
	d \tau_P(A) = T^* \, A \, T, \quad
		d \tau_P(B) = T^* \, B \, T,
\end{equation}
and they are viewed as tangent vectors to $\mathcal{S}_+({\bf R}^n)$
or $\mathcal{S}_+({\bf C}^n)$ at $\tau(P)$.  One can check that
\begin{equation}
	\langle d \tau_P(A), d \tau_P(B) \rangle_{\tau(P)}
		= \langle A, B \rangle_P,
\end{equation}
which is to say that the Riemannian metrics on $\mathcal{S}_+({\bf
R}^n)$, $\mathcal{S}_+({\bf C}^n)$ are invariant under the actions
of $GL({\bf R}^n)$, $GL({\bf C}^n)$ that we have defined.

	Of course the usual flat Riemannian metrics are defined
as follows.  If $T$ is an element of $\mathcal{S}({\bf R}^n)$
or of $\mathcal{S}({\bf C}^n)$, then again two tangent vectors
$A$, $B$ to $\mathcal{S}({\bf R}^n)$ or $\mathcal{S}({\bf C}^n)$
at $T$ are given by self-adjoint linear transformations on
${\bf R}^n$ or ${\bf C}^n$, as appropriate, and their inner
product in the standard flat Riemannian metric is given by
the usual inner product, namely
\begin{equation}
	\tr A \, B.
\end{equation}
The flatness of this Riemannian metric is reflected in the fact
that it does not depend on the point $T$ in the space.

	The metrics that we have defined on $\mathcal{S}_+({\bf
R}^n)$, $\mathcal{S}_+({\bf C}^n)$ reduce exactly to the standard flat
metric at the identity operator $I$.  The differential of the
exponential function at the zero operator is equal to the identity
mapping, 
\begin{equation}
	d \exp_0 (A) = A,
\end{equation}
as we have noted previously, and thus the standard flat metric on the
tangent space of $\mathcal{S}({\bf R}^n)$ or $\mathcal{S}({\bf C}^n)$
at $0$ agrees with the metric that we have defined on
$\mathcal{S}_+({\bf R}^n)$ or $\mathcal{S}_+({\bf C}^n)$, respectively,
at the identity operator $I$ with respect to the correspondence
given by the differential of the exponential function.  We would
like to show that in fact this agreement works to another term
in the Taylor expansion.

	Basically this means that if $T$ is a self-adjoint linear
transformation on ${\bf R}^n$ or on ${\bf C}^n$, which we think of as
being near the origin, and if $A$, $B$ are two self-adjoint linear
transformations on the same space, which we think of as tangent
vectors to $\mathcal{S}({\bf R}^n)$ or $\mathcal{S}({\bf C}^n)$ at
$T$, then the images of $A$, $B$ under the differential of the
exponential map at $T$ have inner product with respect to the
Riemannian metric defined above at the exponential of $T$ which agrees
to second order with the inner product of $A$ in the standard flat
metric on $\mathcal{S}({\bf R}^n)$ or $\mathcal{S}({\bf C}^n)$.
Explicitly, this means that
\begin{eqnarray}
\lefteqn{\langle d \exp_T(A), d \exp_T(B) \rangle_{\exp T}}	\\
	   & & 	= \tr (\exp (-T)) (d \exp_T(A)) (\exp (-T)) (d \exp_T (B))
\end{eqnarray}
agrees with
\begin{equation}
	\tr A \, B
\end{equation}
up to terms of order $O(\|T\|^2)$.  This is not difficult to check,
using the facts that $\exp (-T) = I - T +O(\|T\|^2)$,
\begin{equation}
	d \exp_T (A) = A + \frac{1}{2}(T \, A + A \, T) + O(\|T\|^2),
\end{equation}
and similarly for $B$.

	As a consequence, if $C$ is any self-adjoint linear transformation
on ${\bf R}^n$ or ${\bf C}^n$, then
\begin{equation}
	\exp (t \, C)
\end{equation}
satisfies the equation for geodesics in the space of positive-definite
linear transformations at $t = 0$.  Basically this is because the line
$t \, C$ satisfies the equation for geodesics in the flat space of
self-adjoint linear transformations at the origin.  Of course we are
also using the fact that the standard flat metric on the space of
self-adjoint linear transformations around $0$ agrees with the
Riemannian metric on the space of positive-definite metrics around the
identity operator $I$ with respect to the exponential mapping as well 
as they do.

	In fact, we get that
\begin{equation}
	\exp (t \, C)
\end{equation}
satisfies the equation for geodesics in the space of positive-definite
matrices for all real numbers $t$.  For a fixed real number $t_0$, the
statement that this curve satisfies the equation for geodesics at
$t_0$ is equivalent to the statement that $\exp ((t_0 + t) \, C)$
satisfies the equation for geodesics at $t = 0$.  This is in turn
equivalent to the statement that $T^* \, \exp (t \, C) \, T$ satisfies
the equation for geodesics at $t = 0$ when $T = \exp (t_0 \, C / 2)$,
and this follows from the fact that $\exp (t \, C)$ satisfies the
equation for geodesics at $t = 0$, since transformations of the form
$P \mapsto T^* \, P \, T$ are isometries on the spaces of
positive-definite matrices and therefore preserve geodesics.

	Thus exponentials of straight lines through the origin in the
spaces of self-adjoint linear transformations give rise to geodesics
through the identity operator $I$ in the corresponding spaces of
positive-definite linear transformations.  By standard results about
uniqueness of initial value problems for ordinary differential
equations this accounts for all of the geodesics in the spaces of
positive definite linear transformations which pass through the
identity operator $I$.  Using mappings of the form $P \mapsto T^* \, P
\, T$ one can move these curves to other places and thereby account
for all of the geodesics in the spaces of positive-definite linear
transformations.

\end{document}